\documentclass{amsart}
\usepackage{amsfonts}
\usepackage{amssymb}
\usepackage{amsbsy}

\newtheorem{theorem}{Theorem}[section]

\newtheorem{corollary}[theorem]{Corollary}
\newtheorem{definition}[theorem]{Definition}

\newtheorem{lemma}[theorem]{Lemma}

\newtheorem{remark}[theorem]{Remark}

\numberwithin{equation}{section}
\begin{document}
\title[A note on $LP$-Kenmotsu manifolds admitting Ricci-Yamabe solitons]{A note on $LP$-Kenmotsu manifolds admitting Ricci-Yamabe solitons} 
\author{ Mobin Ahmad, Gazala and Mohd Bilal}
\maketitle

\begin{abstract} In the current note, we study Lorentzian para-Kenmotsu (in brief, $LP$-Kenmotsu) manifolds admitting Ricci-Yamabe solitons  (RYS) and gradient Ricci-Yamabe soliton (gradient RYS).  At last by constructing a 5-dimensional non-trivial example we illustrate our result.
\end{abstract}

 \medskip


\section {\bf Introduction}

%In response to his own work on Ricci flow, Hamilton introduced the notion of Yamabe flow \cite{ HRS1,  HRS2} which deforms the %metric of a Riemannian manifold through the governing equation 
%\begin{eqnarray*}
%\frac {\partial}{\partial t}g(t)=-r (t) g(t), \ \quad g(0)=g_0,
%\end{eqnarray*}
%where $r(t)$ is the scalar curvature of the metric $g(t)$. For the manifold of dimension $2$, the Yamabe flow is equivalent to the %Ricci flow defined by $\frac {\partial}{\partial t}g(t)=-2 S(g(t)),$ where $S$ denotes the Ricci tensor. However, for the case of %higher dimension ($>2$), the Yamabe and Ricci flows do not coincide, since the Yamabe flow preserves the conformal class of %$g(t)$ but the Ricci flow does not in general. The Ricci solitons and Yamabe solitons correspond to self similar solutions of the Ricci %flow and Yamabe flow are given by \cite{BR1, BR2}
%\begin{eqnarray*}
%\pounds_V g+2S+2 \lambda g=0,
%\end{eqnarray*}
%and
%\begin{eqnarray*}
%\pounds_V g=2(r-\lambda) g,
%\end{eqnarray*}
%respectively; where $\pounds_V$ represents the Lie derivative operator along the smooth vector field $V$ (called soliton vector %field) on $M$, $\lambda \in R$ (called soliton constant of the manifold $M$) and $S$ is the Ricci tensor.

In 2019, a scalar combination of Ricci and Yamabe flows was proposed by the authors G{\"u}ler and Crasmareanu \cite{GC}, this advanced class of geometric flows called Ricci-Yamabe  (RY) flow of type $(\sigma, \rho)$ and is defined by 
\begin{eqnarray*}
\frac {\partial}{\partial t}g(t)+2 \sigma S(g(t))+\rho r(t) g(t)=0,  \ \quad g(0)=g_0
\end{eqnarray*}
for some scalars $\sigma$ and $\rho$. 

A solution to the RY flow is called  RYS if it depends
only on one parameter group of diffeomorphism and scaling. A Riemannian (or semi-Riemannian) manifold $M$ is said to have a RYS if 
\begin{eqnarray} \label{1.1}
\pounds_K g+2 \sigma S+(2 \Lambda-\rho r) g=0,
\end{eqnarray}
where $ \sigma, \rho, \Lambda \in \mathbb R$ (the set of real numbers). If $K$ is the gradient of a smooth function $v$ on $M$, then \eqref{1.1} is called the gradient Ricci-Yamabe soliton (gradient RYS) and hence  \eqref{1.1} turns to
\begin{eqnarray} \label{1.2}
\nabla^2v+\sigma S+( \Lambda-\frac{\rho r}{2})g=0,
\end{eqnarray}
where $\nabla^2v$ is the Hessian of $v$. It is to be noted that a RYS of types $(\sigma,0)$ and $(0, \rho)$ are known as $\sigma-$Ricci soliton and $\rho-$Yamabe soliton, respectively. A RYS is said to be shrinking , steady  or expanding  if $\Lambda<0, =0$ or $>0$, respectively.  A RYS is said to
be a \\
$\bullet$ Ricci soliton \cite{HRS1} if $\sigma=1, \rho=0$,\\
$\bullet$  Yamabe soliton \cite{HRS2} if $\sigma=0, \rho=1$,\\
$\bullet$  Einstein soliton \cite{CM1} if $\sigma=1, \rho=-1$,\\

As a continuation of this study, we tried to study RYS in the frame-work of $LP$-Kenmotsu manifolds of dimension $n$. We recommend the papers \cite{BAM1, BAM2, CV, HD, HA,  HPM, LH, PCP, SK, VK, YHI} and the references therein for more details about the related studies.

\section{\bf Preliminaries}
An $n$-dimensional differentiable manifold $M$ with structure $(\varphi, \zeta,\nu,g)$ is said to be a  Lorentzian almost paracontact metric manifold, if it admits a $(1,1)$-tensor field $\varphi$, a contravariant vector field $\zeta$, a 1-form $\nu$ and a Lorentzian metric $g$ satisfying
\begin{equation}\label{2.1}
{\nu {(\zeta )}+1=0},
\end{equation}
\begin{equation}\label{2.2}
{\varphi }^{2}{E}={E}+{\nu }(E){\zeta},
\end{equation}
 \begin{equation}\label{2.3}
\varphi \zeta=0, \quad \nu(\varphi E)=0,
\end{equation}
\begin{equation}\label{2.4}
g(\varphi {E},\varphi {F})=g(E,F)+ {\nu }(E){\nu }(F),
\end{equation}
\begin{equation}\label{2.5}
  g(E,{\zeta })={\nu }(E),
 \end{equation}
\begin{equation}\label{2.6}
 \varphi (E,F)=\varphi (F,E)=g(E,\varphi F)
 \end{equation}
for any vector fields $E,F\in \chi(M)$, where $\chi(M)$ is the Lie algebra of vector fields on $M$.\\
If $\zeta$ is a killing vector field, the (para) contact structure is called a $K$-(para) contact. In such a case, we have
 \begin{equation}\label{2.7}
\nabla_{E}\zeta=\varphi E.
\end{equation}
Recently, the authors Haseeb and Prasad defined and studied the following notion: 
\begin{definition}
A Lorentzian almost paracontact manifold $M$ is called
Lorentzian para-Kenmostu manifold if \cite{HP1} 
\begin{equation}\label{2.9}
(\nabla _{E}\varphi )F=-g(\varphi E,F)\zeta -\nu (F)\varphi E
\end{equation}
for any $E,F$ on $M.$
\end{definition}
\noindent
In an $LP$-Kenmostu manifold, we have
\begin{equation}\label{2.10}
\nabla _{E}\zeta =-E-\nu (E)\zeta,
\end{equation}
\begin{equation}\label{2.11}
(\nabla _{E}\nu )F=-g(E,F)-\nu (E)\nu (F),
\end{equation}%
where $\nabla $ denotes the Levi-Civita connection respecting to the Lorentzian metric $g$.

\noindent
Furthermore, in an $LP$-Kenmotsu manifold, the following relations hold \cite{HP1}:
\begin{equation}\label{2.12}
g(R(E,F)G,\zeta)=\nu(R(E,F)G)=g(F,G)\nu(E)-g(E,G)\nu(F),
\end{equation}
\begin{equation}\label{2.13}
R(\zeta,E)F=-R(E,\zeta)F=g(E,F)\zeta-\nu(F)E,
\end{equation}
\begin{equation}\label{2.14}
R(E,F)\zeta=\nu(F)E-\nu(E)F,
\end{equation}
\begin{equation}\label{2.15}
R(\zeta,E)\zeta=E+\nu(E)\zeta,
\end{equation}
\begin{equation}\label{2.16}
S(E,\zeta)=(n-1)\nu(E),~~S(\zeta,\zeta)=-(n-1),
\end{equation}
\begin{equation}\label{2.17}
Q \zeta=(n-1) \zeta
\end{equation}
 for any $E,F,G \in\chi(M)$, where $R, S$ and $Q$ represent the curvature tensor, the Ricci tensor and the $Q$ Ricci operator, respectively.
\begin{definition} \cite{YK2}
An $LP$-Kenmotsu manifold $M$ is said to be $\nu$-Einstein manifold if its $S (\neq 0)$ is of the form
\begin{equation}\label{2.18}
S(E,F)=ag(E,F)+b\nu(E)\nu(F),
\end{equation}
where $a$ and $b$ are smooth functions on $M$. In particular, if $b = 0$, then $M$ is termed as an Einstein manifold.
\end{definition} 

\noindent
\begin{remark} \cite{HP2} In an $LP$-Kenmotsu manifold of $n$-dimension, $S$ is of the form
\begin{equation} \label{2.19}
S(E,F)=(\frac {r}{n-1}-1)g(E,F)+(\frac {r}{n-1}-n)\nu(E)\nu(F),
\end{equation}
where $r$ is the scalar curvature of the manifold.
\end{remark}
\begin{lemma} In an $n$-dimensional $LP$-Kenmotsu manifold, we have
\begin{equation} \label{2.20}
\zeta(r)=2(r-n(n-1)),
\end{equation}
\begin{equation} \label{2.21}
(\nabla_E Q)\zeta=QE-(n-1)E,
\end{equation}
\begin{equation} \label{2.22}
(\nabla_{\zeta} Q)E=2QE-2(n-1)E
\end{equation}
for any $E$ on $M$.
\begin{proof} Equation \eqref{2.19} yields
\begin{equation} \label{2.23}
QE=(\frac {r}{n-1}-1)E+(\frac {r}{n-1}-n)\nu(E)\zeta.
\end{equation}
Taking the covariant derivative of \eqref{2.23} with respect to $F$ and making use of \eqref{2.10} and \eqref{2.11}, we lead to
\begin{eqnarray*} 
(\nabla_FQ)E=\frac{F(r)}{n-1}(E+\nu(E)\zeta)
-(\frac {r}{n-1}-n)(g(E,F)\zeta+\nu(E)F+2 \nu(E)\nu(F)\zeta).
\end{eqnarray*}
By contracting $F$ in the foregoing equation and using trace $\{F  \to (\nabla_F Q)E\}=\frac{1}{2}E(r)$, we find
\begin{eqnarray*} \label{2.24}
\frac {n-3}{2(n-1)} E(r)=\big\{\frac{\zeta(r)}{n-1}-(r-n(n-1))\big\} \nu(E),
\end{eqnarray*}
which by replacing $E$ by $\zeta$ and using \eqref{2.1} gives \eqref{2.20}. We refer the readers to see \cite {LHM} for the proof of \eqref{2.21} and \eqref{2.22}.
\end{proof}
\end{lemma}
\begin{remark}
From the equation $(\ref{2.20})$, it is noticed that if an $n$-dimensional $LP$-Kenmotsu manifold possesses the constant scalar curvature, then  $r=n(n-1)$ and hence \eqref{2.19} reduces to $S(E,F)=(n-1)g(E,F)$. Thus, the manifold under consideration
is an Einstein manifold. 
\end{remark}

\section{\bf Ricci-Yamabe solitons on $LP$-Kenmotsu manifolds}
 Let the metric of an $n$-dimensional $LP$-Kenmotsu manifold be a Ricci-Yamabe soliton $(g,K,\Lambda, \sigma, \rho)$, then $(\ref{1.1})$ holds. By differentiating $(\ref{1.1})$ covariantly with resprct to $G$, we have
\begin{eqnarray}\label{3.1}
(\nabla_G\pounds_K g)(E,F)&=&-2 \sigma (\nabla_G S)(E,F)+\rho (Gr) g(E,F).
\end{eqnarray}  
 Since $\nabla g=0$, then the following formula \cite{YK1}
\begin{equation*}
(\pounds_K \nabla_E g-\nabla_E \pounds_K  g-\nabla_{[K,E]} g)(F,G)=-g((\pounds_K \nabla)(E,F),G)-g((\pounds_K \nabla)(E,G),F)
\end{equation*}
turns to
\begin{equation*}
(\nabla_E \pounds_K g)(F,G)=g((\pounds_K \nabla)(E,F),G)+g((\pounds_K \nabla)(E,G),F).
\end{equation*}
Since the operator $\pounds_K \nabla$ is symmetric, therefore we have
\begin{equation*}
2g((\pounds_K \nabla)(E,F),G)=(\nabla_E \pounds_K g)(F,G)+(\nabla_F \pounds_K g)(E,G)-(\nabla_G \pounds_K g)(E,F),
\end{equation*}
which by using \eqref{3.1} takes the form
\begin{eqnarray}\label{3.2}
 2g((\pounds_K \nabla)(E,F),G)&=&-2 \sigma [(\nabla_E S)(F,G)+ (\nabla_F S)(G,E)+(\nabla_G S)(E,F)] \notag  \\ 
 &&+\rho[(E r) g(F,G)+ (F r) g(G,E)+ (G r) g(E,F)].
\end{eqnarray}
Putting $F=\zeta$ in \eqref{3.2} and using \eqref{2.5}, we find
\begin{eqnarray}\label{3.3}
2g((\pounds_K \nabla)(E,\zeta),G)&=&-2 \sigma[ (\nabla_E S)(\zeta,G)+(\nabla_{\zeta} S)(G,E)-(\nabla_G S)(E,\zeta)]\notag\\
&&+\rho[ (E r) \nu(G)+2(r-n(n-1))g(E,G)-(Gr)\nu(E)]
\end{eqnarray}
By virtue of \eqref{2.21} and \eqref{2.22}, \eqref{3.3} leads to
\begin{eqnarray*}
2g((\pounds_K \nabla)(E,\zeta),G)&=&-4 \sigma [S(E,G)-(n-1) g(E,G)]\notag  \\
&&+\rho[ (E r) \nu(G)+2(r-n(n-1))g(E,G)-(Gr)\nu(E)].
\end{eqnarray*}
By eliminating $G$ from the foregoing equation, we have
\begin{eqnarray}\label{3.4}
2 (\pounds_K \nabla)(F,\zeta)&=&\rho g(Dr,F) \zeta-\rho (Dr) \nu (F)-4 \sigma QF\\
&&+[4 \sigma (n-1)+2 \rho (r-n(n-1))]F. \notag
\end{eqnarray}
If we take $r$ as constant, then from \eqref{2.20} we find $r=n(n-1)$, and hence \eqref{3.4} reduces to
\begin{eqnarray}\label{3.5}
(\pounds_K \nabla)(F,\zeta)&=&-2 \sigma QF+2 \sigma (n-1) F.
\end{eqnarray}
Taking covariant derivative of \eqref{3.5} with respect to $E$, we have
\begin{eqnarray}\label{3.6}
(\nabla_E \pounds_K \nabla)(F,\zeta)&=&(\pounds_K \nabla)(F,E)-2 \sigma \nu(E)[QF-(n-1)F]\\
&-&2 \sigma (\nabla_E Q) F. \notag
\end{eqnarray}
 Again from \cite{YK1}, we have
\begin{eqnarray*}
 (\pounds_K R)(E,F)G=(\nabla_E\pounds_K \nabla)(F,G)-(\nabla_F\pounds_K \nabla)(E,G),
\end{eqnarray*}
which by putting $G=\zeta$ and using \eqref{3.6} takes the form
\begin{eqnarray}\label{3.7}
 (\pounds_K R)(E,F)\zeta&=&2 \sigma \nu(F)(QE-(n-1)E)-2 \sigma \nu(E)(QF-(n-1)F)\\
&&-2 \sigma ((\nabla_E Q) F- (\nabla_F Q) E). \notag
\end{eqnarray}
Putting $F=\zeta$ in \eqref{3.7} then using \eqref{2.1}, \eqref{2.2}, \eqref{2.21} and \eqref{2.22}, we arrive at
\begin{eqnarray}\label{3.8}
 (\pounds_K R)(E,\zeta)\zeta = 0.
\end{eqnarray}
The Lie derivative of $R(E,\zeta)\zeta=-E-\nu(E)\zeta$ along $K$ leads to
\begin{eqnarray}\label{3.9}
 (\pounds_K R)(E,\zeta)\zeta-g(E,\pounds_K \zeta)\zeta+2 \nu(\pounds_K \zeta)E =-(\pounds_K \nu) (E)\zeta.
\end{eqnarray}
From \eqref{3.8} and \eqref{3.9}, we have
\begin{eqnarray}\label{3.10}
 (\pounds_K \nu) (E)\zeta=-2 \nu(\pounds_K \zeta)E+g(E,\pounds_K \zeta)\zeta.
\end{eqnarray}
Taking the Lie derivative of $g(E,\zeta)=\nu(E)$, we find
\begin{eqnarray}\label{3.11}
 (\pounds_K \nu) (E)=g(E,\pounds_K \zeta)+ (\pounds_K g)(E,\zeta).
\end{eqnarray}
By putting $F=\zeta$ in \eqref{1.1} and using \eqref{2.16}, we have
\begin{eqnarray}\label{3.12}
(\pounds_K g)(E,\zeta)=-\{2\sigma (n-1)+2 \Lambda-\rho n(n-1)\}\nu(E),
\end{eqnarray}
where $r=n(n-1)$ being used.\\
The Lie derivative of $g(\zeta,\zeta)=-1$ along $K$ we lead to
\begin{eqnarray}\label{3.13}
(\pounds_K g)(\zeta,\zeta)=-2 \nu (\pounds_K \zeta).
\end{eqnarray}
From \eqref{3.12} and \eqref{3.13}, we find
\begin{eqnarray}\label{3.14}
\nu (\pounds_K \zeta)=-\{\sigma (n-1)+\Lambda-\frac{\rho n(n-1)}{2}\}.
\end{eqnarray}
Now, combining the equations \eqref{3.10}, \eqref{3.11},  \eqref{3.12} and \eqref{3.14}, we find
\begin{eqnarray}\label{3.15}
\Lambda=\frac{\rho n(n-1)}{2}-\sigma (n-1).
\end{eqnarray}
Thus, we have
\begin{theorem}
 Let $(M,g)$ be an $n$-dimensional $LP$-Kenmotsu manifold admitting Ricci-Yamabe soliton $(g,K,\Lambda,\sigma, \rho)$ with constant scalar curvature tensor, then $\Lambda=\frac{\rho n(n-1)}{2}-\sigma (n-1).$
\end{theorem}

For $\sigma=1$ and $\rho=0$, from \eqref{3.15} we have $\Lambda=-(n-1)$. Thus, we have the following:
\begin{corollary}
 If an $n$-dimensional $LP$-Kenmotsu manifold admits a Ricci soliton with constant scalar curvature, then the soliton is  shrinking.
\end{corollary}

For $\sigma=0$ and $\rho=1$, from \eqref{3.15} we have $\Lambda=\frac{n(n-1)}{2}$. Thus, we have the following:
\begin{corollary}
 If an $n$-dimensional $LP$-Kenmotsu manifold admits a Yamabe soliton with constant scalar curvature, then the soliton is shrinking.
\end{corollary}

For $\sigma=1$ and $\rho=-1$, from \eqref{3.15} we have $\Lambda=-\frac{(n^2-1)}{2}$. Thus, we have the following:
\begin{corollary}
 If an $n$-dimensional $LP$-Kenmotsu manifold admits an Einstein soliton with constant scalar curvature, then the soliton is shrinking.
\end{corollary}

Now, we consider the metric of an $n$-dimensional $LP$-Kenmotsu manifold as a Ricci-Yamabe soliton $(g,\zeta,\Lambda, \sigma, \rho)$, then from $(\ref{1.1})$ and \eqref{2.10} we have
\begin{eqnarray} \label{3.13}
 \ \ \ \ S(E,F)=-\frac{1}{\sigma} (\Lambda-1-\frac{\rho r}{2})g(E,F)+\frac{1}{\sigma} \nu(E)\nu(F), \ \ where  \ \ \sigma \neq0.
\end{eqnarray}
By putting $F=\zeta$ in \eqref{3.13} and using \eqref{2.16}, we find
\begin{eqnarray} \label{3.14}
\Lambda=\frac{\rho r}{2}-\sigma (n-1).
\end{eqnarray}
Now, comparing \eqref{2.19} and \eqref{3.14}, we have $r=\frac{n-1}{\sigma}+n(n-1)$, which by using in \eqref{3.14} it follows that $\Lambda=-\sigma (n-1)+\frac{\rho  (n-1)(1+n \sigma)}{2 \sigma}.$
Thus, we have the following theorem:
\begin{theorem}
An $n$-dimensional $LP$-Kenmotsu manifold with constant scalar curvature admitting Ricci-Yamabe soliton $(g,\zeta,\Lambda,\sigma, \rho)$ is an $\nu$-Einstein manifold. Moreover, the soliton is expanding, steady or shrinking according to $\frac {\rho}{\sigma}> 2 \sigma -\rho n$,  $\frac {\rho}{\sigma}=2 \sigma -\rho n$, or  $\frac {\rho}{\sigma}< 2 \sigma -\rho n$.
\end{theorem}

\section{\protect\vspace{0.2cm}\textbf{Gradient Ricci-Yamabe solitons on $LP$-Kenmotsu manifolds\label{sect-Curvature}}}
\begin{definition} A Riemannian $($or semi-Riemannian$)$ metric $g$ on $M$ is called a gradient RYS, if
\begin{eqnarray} \label{4.1}
Hess v+\sigma S+(\Lambda-\frac{\rho r}{2}) g=0,
\end{eqnarray}
where $Hess v$ denotes the Hessian of a smooth function $v$ on $M$ and defined by $Hess v=\nabla \nabla v$.
\end{definition}

Let $M$ be an $n$-dimensional $LP$-Kenmotsu manifold with $g$ as a gradient RYS. Then
equation (4.1) can be written as
\begin{eqnarray}\label{4.2}
\nabla_E D v+\sigma Q E+(\Lambda-\frac {\rho r}{2})E=0,
\end{eqnarray}
for all vector fields $E$ on $M$, where $D$ denotes the gradient operator of $g$. 
Taking the covariant derivative of \eqref{4.2} with respect to $F$, we have
\begin{eqnarray}\label{4.3}
\ \  \ \nabla_F\nabla_E D v=-\sigma\{(\nabla_F Q) E+Q(\nabla_F E)\} +\rho \frac{F(r)}{2}E-(\Lambda-\frac {\rho r}{2})\nabla_F E.
\end{eqnarray}
Interchanging $E$ and $F$ in \eqref{4.3}, we lead to
\begin{eqnarray}\label{4.4}
\ \  \ \nabla_E\nabla_F D v=-\sigma\{(\nabla_E Q) F+Q(\nabla_E F)\} +\rho \frac{E(r)}{2}F-(\Lambda-\frac {\rho r}{2})\nabla_E F.
\end{eqnarray}
By making use of \eqref{4.2}-\eqref{4.4}, we find
\begin{eqnarray}\label{4.5}
R(E,F) D v=\sigma\{(\nabla_F Q) E-(\nabla_E Q) F\} +\frac{\rho}{2} \{ E(r) F-F(r)E\}.
\end{eqnarray}
Now, from \eqref{2.19}, we find
\begin{eqnarray*}
QE=(\frac {r}{n-1}-1)E+(\frac {r}{n-1}-n)\nu(E)\zeta,
\end{eqnarray*}
which on taking covariant derivative with repect to $F$ leads to
\begin{eqnarray}\label{4.6}
(\nabla_F Q) E&=&\frac {F(r)}{n-1}(E+\nu(E)\zeta)-(\frac {r}{n-1}-n)(g(E,F)\zeta\\
&&+2 \nu(E)\nu(F)\zeta+\nu(E)F). \notag
\end{eqnarray}
By using \eqref{4.6} in \eqref{4.5}, we have
\begin{eqnarray}\label{4.7}
R(E,F) D v&=&\frac{(n-1)\rho-2 \sigma}{2(n-1)} \{E(r) F-F(r)E\}+\frac{\sigma}{n-1}\{F(r)\nu(E)\zeta-E(r)\nu(F)\zeta\} \notag \\
&&-\sigma(\frac{r}{n-1}-n)(\nu(E)F-\nu(F)E).
\end{eqnarray}
Contracting forgoing equation along $E$ gives
\begin{eqnarray}\label{4.8}
S(F, D v)&=&\big\{\frac {(n-1)^2\rho-2\sigma (n-2)}{n-1}\big \} F(r)\\
&&+\frac{\sigma(n-3) (r-n(n-1))}{n-1}\nu(F). \nonumber
\end{eqnarray}
From the equation \eqref{2.19}, we can write
\begin{eqnarray}\label{4.9}
S(F, D v)=(\frac {r}{n-1}-1) F(v)+(\frac {r}{n-1}-n)\nu(F) \zeta (v).
\end{eqnarray}
Now, by equating \eqref{4.8} and \eqref{4.9}, then putting $F=\zeta$ and using \eqref{2.1}, \eqref{2.20}, we find
\begin{eqnarray}\label{4.10}
 \zeta (v)= \frac {r-n(n-1)}{n-1}\{2 (n-1)\rho-\frac{\sigma(5n-13)}{n-1}\}.
\end{eqnarray}
Taking the inner product of  \eqref{4.7} with $\zeta$, we get
\begin{eqnarray*}
 F(v) \nu(E)-E(v)\nu(F)= \frac {\rho}{2}\{E(r)\nu(F)-F(r)\nu(E)\},
\end{eqnarray*}
which by replacing $E$ by $\zeta$ and using \eqref{2.20}, \eqref{4.10}, we infer
\begin{eqnarray}\label{4.11}
 F(v)=- (r-n(n-1)) \{3 \rho-\frac{\sigma(5n-13)}{(n-1)^2}\}\nu(F)-\frac {\rho}{2}F(r).
\end{eqnarray}
If we take $r$ as constant, then from Remark 2.5, we get $r=n(n-1)$. Thus, \eqref{4.11} leads to $F(v)=0$. This implies that $v$ is constant. Thus, the soliton under the consideration is trivial. Hence we state:
\begin{theorem} \label{thm4.1} 
If the metric of an $LP$-Kenmotsu manifold of constant scalar curvature tensor admitting a special type of vector field is gradient RYS, then the soliton is trivial.
\end{theorem}
For $v$ constant, \eqref{1.2} turns to
\begin{eqnarray*}
 \sigma QE=-(\Lambda-\frac{\rho r}{2})E,
\end{eqnarray*}
which leads to
\begin{eqnarray}\label{4.12}
 S(E,F)=-\frac{1}{\sigma} (\Lambda-\frac{\rho n(n-1)}{2})g(E,F), \ \ \ \sigma \neq 0.
\end{eqnarray}
By putting $E=F=\zeta$ in \eqref{4.12} and using \eqref{2.16}, we obtain
 \begin{eqnarray}\label{4.13}
\Lambda=\frac{\rho n (n-1)}{2}-\sigma(n-1).
\end{eqnarray}
\begin{corollary}
 If an $n$-dimensional $LP$-Kenmotsu manifold admits a gradient Ricci soliton with the constant scalar curvature, then the manifold under the consideration is an Einstein manifold and $\Lambda=\frac{\rho n (n-1)}{2}-\sigma(n-1).$
\end{corollary}
For $\sigma=1$ and $\rho=0$, from \eqref{4.13} we find $\Lambda=-(n-1)$. Thu, we have the following:
\begin{corollary}
 If an $n$-dimensional $LP$-Kenmotsu manifold admits a gradient Ricci soliton with the constant scalar curvature, then the soliton is  shrinking.
\end{corollary}
For $\sigma=1$ and $\rho=-1$, from \eqref{4.13} we have $\Lambda=-\frac{(n-1)(n+2)}{2}$. Thus, we have the following:
\begin{corollary}
 If an $n$-dimensional $LP$-Kenmotsu manifold admits an gradient Einstein soliton with constant scalar curvature, then the soliton is shrinking.
\end{corollary}

\noindent
{\bf Example.} We consider the 5-dimensional manifold $M^5=\left\{(x_1,x_2,x_3,x_4,x_5)\in \mathbb R^5: x_5>0\right\}$, where $(x_1,x_2,x_3,x_4,x_5)$ are the standard coordinates in $\mathbb R^5$. Let $\varrho_1,$ $\varrho_2$, $\varrho_3$, $\varrho_4$ and $\varrho_5$ be the vector fields on $M^5$ given by
\begin{equation*}
\varrho_1=e^{x_5}\frac{\partial}{\partial x_1}, ~~\varrho_2=e^{x_5}\frac{\partial}{\partial {x_2}}, ~~\varrho_3=e^{x_5}\frac{\partial}{\partial {x_3}},~~\varrho_4=e^{x_5}\frac{\partial}{\partial x_4},~~\varrho_5=\frac{\partial}{\partial x_5}=\zeta,
\end{equation*}
which are linearly independent at each point of $M^5$. Let $g$ be the Lorentzian metric defined by
\begin{equation*}
\ \ \ \ \ \ \ \ \ \ \ \ \ \ \ \  g(\varrho_i, \varrho_i)=1,\ \quad {\text {for}}\ \quad 1\leq i\leq 4\ \quad {\text {and}}  \ \quad g(\varrho_5, \varrho_5)=-1,
\end{equation*}
\begin{equation*}
g(\varrho_i, \varrho_j)=0, \ \quad {\text {for}}\ \quad  i \neq  j, \ \quad 1\leq i,j\leq 5.
\end{equation*}
Let $\nu$ be the 1-form defined by $\nu(E)=g(E,\varrho_5)=g(\varrho,\zeta)$ for all $E\in \chi {(M^5)}$,
and let $\varphi$ be the $(1,1)$-tensor field defined by
\begin{equation*}
\varphi \varrho_1=-\varrho_2, ~~\varphi \varrho_2=-\varrho_1, ~~\varphi\varrho_3=-\varrho_4, ~~\varphi \varrho_4=-\varrho_3,~~\varphi \varrho_5=0.
\end{equation*}
 By applying linearity of $\varphi$ and $g$, we have
\begin{equation*}
\nu(\zeta)=g(\zeta,\zeta)=-1, ~~\varphi^2 E=E+\nu(E)\zeta~~{\text {and}} ~~g(\varphi E, \varphi F)=g(E,F)+\nu(E)\nu(F)
\end{equation*}
for all $E, F\in \chi{(M^5)}$. Thus for $\varrho_5=\zeta$, the structure $(\varphi, \zeta,\nu,g)$ defines a Lorentzian almost paracontact metric structure on $M^5$.  Then we have
\begin{equation*}
[\varrho_i,\varrho_j]=-\varrho_i, \ \quad {\text {for}}\ \quad 1\leq i\leq 4, j=5,   \ \ \ 
\end{equation*}
\begin{equation*}
 [\varrho_i,\varrho_j]=0, \ \quad {\text {otherwise}}.\ \ \ \ \ \ \ \ \ \ \ \ \ \ \ \ \ \ \ \ \ 
\end{equation*}

By using  Koszul's formula, we can easily find
we obtain
\begin{align*} \nabla_{\varrho_i} {\varrho}_j=
	\begin{cases}
                     -\varrho_5,\ \quad   1 \le i=  j\le 4,\\
		-\varrho_i,\ \quad   1 \le i\le 4,  j= 5,\\
		0 ,\ \quad otherwise.
	\end{cases}
\end{align*}
Also one can easily verify that 
\begin{equation*}
\nabla_E \zeta=-E-\eta(E)\zeta  \ \quad {\text {and}}\ \quad (\nabla_E \varphi)F=-g(\varphi E,F)\zeta-\nu(F)\varphi E.
\end{equation*}
Therefore, the manifold is an $LP$-Kenmotsu manifold.\\
\noindent
From the above results, we can easily obtain the non-vanishing components of  $R$ as follows:
\begin{equation*}
R(\varrho_1,\varrho_2)\varrho_1=-\varrho_2,~~~R(\varrho_1,\varrho_2)\varrho_2=\varrho_1,~~~R(\varrho_1,\varrho_3)\varrho_1=-\varrho_3,~~~ R(\varrho_1,\varrho_3)\varrho_3=\varrho_1,
\end{equation*}
\begin{equation*}
 R(\varrho_1,\varrho_4)\varrho_1=-v_4,~~~R(\varrho_1,\varrho_4)\varrho_4=\varrho_1,~~~R(\varrho_1,\varrho_5)\varrho_1=-\varrho_5,~~~ R(\varrho_1,\varrho_5)\varrho_5=-\varrho_1,
\end{equation*}
\begin{equation*}
 R(\varrho_2,\varrho_3)\varrho_2=-\varrho_3,~~~R(\varrho_2,\varrho_3)\varrho_3=\varrho_2,~~~R(\varrho_2,\varrho_4)\varrho_2=-\varrho_4,~~~ R(\varrho_2,\varrho_4)\varrho_4=\varrho_2,
\end{equation*}
\begin{equation*}
R(\varrho_2,\varrho_5)\varrho_2=-\varrho_5,~~~R(\varrho_2,\varrho_5)\varrho_5=-\varrho_2,~~~R(\varrho_3,\varrho_4)\varrho_3=-\varrho_4,~~~ R(\varrho_3,\varrho_4)\varrho_4=\varrho_3,
\end{equation*}
\begin{equation*}
 R(\varrho_3,\varrho_5)\varrho_3=-\varrho_5,~~~R(\varrho_3,\varrho_5)\varrho_5=-\varrho_3,~~~R(\varrho_4,\varrho_5)\varrho_4=-\varrho_5,~~~ R(\varrho_4,\varrho_5)\varrho_5=-\varrho_4.
\end{equation*}
Also, we calculate the Ricci tensors as follows:
\begin{equation*}
S(\varrho_1,\varrho_1)=S(\varrho_2,\varrho_2)=S(\varrho_3,\varrho_3)= S(\varrho_4,\varrho_4)=4, \ \ \ \ S(\varrho_5,\varrho_5)=-4.
\end{equation*}
Therefore, we have
 \begin{equation*}
r=S(\varrho_1,\varrho_1)+S(\varrho_2,\varrho_2)+S(\varrho_3,\varrho_3)+S(\varrho_4,\varrho_4)- S(\varrho_5,\varrho_5)=20.
\end{equation*}
Now by taking $Dv=(\varrho_1 v)\varrho_1+(\varrho_2 v)\varrho_2+(\varrho_3 v) \varrho_3+(\varrho_4 v)\varrho_4+( \varrho_5 v) \varrho_5$, we have
\begin{eqnarray*}
\nabla_{\varrho_1} Dv=(\varrho_1(\varrho_1 v)-(\varrho_5 v)){\varrho_1}+{(\varrho_1}(\varrho_2 v)){\varrho_2}+{(\varrho_1}(\varrho_3 v)){\varrho_3}+{(\varrho_1}(\varrho_4 v)){\varrho_4} +(\varrho_1(\varrho_5 v)-(\varrho_1 v)){\varrho_5},
\end{eqnarray*}
\begin{eqnarray*}
\nabla_{\varrho_2} Dv=(\varrho_2(\varrho_1 v)) \varrho_1 +(\varrho_2(\varrho_2 v)-(\varrho_5 v)){\varrho_2}+{(\varrho_2}(\varrho_3 v)){\varrho_3}+{(\varrho_2}(\varrho_4 v)){\varrho_4} +(\varrho_2(\varrho_5 v)-(\varrho_2 v)){\varrho_5},
\end{eqnarray*}
\begin{eqnarray*}
\nabla_{\varrho_3} Dv=(\varrho_3(\varrho_1 v)) \varrho_1 +(\varrho_3(\varrho_2 v)) \varrho_2+(\varrho_3(\varrho_3 v)-(\varrho_5 v)){\varrho_3}+(\varrho_3(\varrho_4 v)){\varrho_4} +(\varrho_3(\varrho_5 v)-(\varrho_3 v)){\varrho_5},
\end{eqnarray*}
\begin{eqnarray*}
\nabla_{\varrho_4} Dv=(\varrho_4(\varrho_1 v)) \varrho_1 +(\varrho_4(\varrho_2 v)) \varrho_2+(\varrho_4(\varrho_3 v)) \varrho_3+(\varrho_4(\varrho_4 v)-(\varrho_5 v)){\varrho_4}+(\varrho_4(\varrho_5 v)-(\varrho_4 v)){\varrho_5},
\end{eqnarray*}
\begin{eqnarray*}
\nabla_{\varrho_5} Dv=(\varrho_5(\varrho_1 v)) \varrho_1 +(\varrho_5(\varrho_2 v)) \varrho_2+(\varrho_5(\varrho_3 v)) \varrho_3+(\varrho_5(\varrho_4 v))\varrho_4+(\varrho_5(\varrho_5 v)){\varrho_5}.
\end{eqnarray*}
Thus, by virtue of \eqref{4.2}, we obtain
\begin{align}\label{4.14}
\begin{cases}
\varrho_1(\varrho_1 v)-\varrho_5 v=-(\Lambda +4 \sigma-10 \rho),\\
\varrho_2(\varrho_2 v)-\varrho_5 v=-(\Lambda +4 \sigma-10 \rho),\\
\varrho_3(\varrho_3 v)-\varrho_5 v=-(\Lambda +4 \sigma-10 \rho),\\
\varrho_4(\varrho_4 v)-\varrho_5 v=-(\Lambda +4 \sigma-10 \rho),\\
\varrho_5(\varrho_5 v)=-(\Lambda +4 \sigma-10 \rho),\\
\varrho_1(\varrho_2 v)=\varrho_1(\varrho_3 v)=\varrho_1(\varrho_4 v)=0,\\
\varrho_2(\varrho_1 v)=\varrho_2(\varrho_3 v)=\varrho_2(\varrho_4 v)=0,\\
\varrho_3(\varrho_1 v)=\varrho_3(\varrho_2 v)=\varrho_3(\varrho_4 v)=0,\\
\varrho_4(\varrho_1 v)=\varrho_4(\varrho_2 v)=\varrho_4(\varrho_3 v)=0,\\
\varrho_1(\varrho_5 v)-(\varrho_1 v)=\varrho_2(\varrho_5 v)-(\varrho_2 v)=0,\\
\varrho_3(\varrho_5 v)-(\varrho_3 v)=\varrho_4(\varrho_5 v)-(\varrho_4 v)=0.
\end{cases}
\end{align}
Thus, the equations in \eqref{4.14} are respectively amounting to
\begin{eqnarray*}
e^{2 x_5}\frac{\partial^2 v}{\partial x_1^2}-\frac{\partial v}{\partial x_5} =-(\Lambda +4 \sigma-10 \rho),
\end{eqnarray*}
\begin{eqnarray*}
e^{2 x_5}\frac{\partial^2 v}{\partial x_2^2}-\frac{\partial v}{\partial x_5} =-(\Lambda +4 \sigma-10 \rho),
\end{eqnarray*}
\begin{eqnarray*}
e^{2 x_5}\frac{\partial^2 v}{\partial x_3^2}-\frac{\partial v}{\partial x_5} =-(\Lambda +4 \sigma-10 \rho),
\end{eqnarray*}
\begin{eqnarray*}
e^{2 x_5}\frac{\partial^2 v}{\partial x_4^2}-\frac{\partial v}{\partial x_5} =-(\Lambda +4 \sigma-10 \rho),
\end{eqnarray*}
\begin{eqnarray*}
\frac{\partial^2 v}{\partial x_5^2} =-(\Lambda +4 \sigma-10 \rho),
\end{eqnarray*}
\begin{eqnarray*}
\frac{\partial^2 v}{\partial x_1 \partial x_2}=\frac{\partial^2 v}{\partial x_1 \partial x_3}=\frac{\partial^2 v}{\partial x_1 \partial x_4} =\frac{\partial^2 v}{\partial x_2 \partial x_3}=\frac{\partial^2 v}{\partial x_2 \partial x_4}=\frac{\partial^2 v}{\partial x_3 \partial x_4} =0,
\end{eqnarray*}
\begin{eqnarray*}
e^{ x_5}\frac{\partial^2 v}{\partial x_5 \partial x_1}+\frac{\partial v}{\partial x_1}=e^{ x_5}\frac{\partial^2 v}{\partial x_5 \partial x_2}+\frac{\partial v}{\partial x_2}=e^{ x_5}\frac{\partial^2 v}{\partial x_5 \partial x_3}+\frac{\partial v}{\partial x_3}=e^{ x_5}\frac{\partial^2 v}{\partial x_5 \partial x_4}+\frac{\partial v}{\partial x_4} =0.
\end{eqnarray*}
From the above  equations it is observed that $v$ is constant for $\Lambda=-4 \sigma+10\rho$. Hence, equation \eqref{4.2} is satisfied. Thus, $g$ is a gradient RYS with the soliton vector field $K=Dv$, where $v$ is constant and $\Lambda=-4 \sigma+10\rho$.
Hence, Theorem \ref{thm4.1} is verified.

\hspace{5mm}

\noindent
Mobin Ahmad\newline
Department of Mathematics, \newline
Integral University, Kursi Road, \newline
Lucknow-226026. \newline
Email : mobinahmad68@gmail.com\newline

\noindent
Gazala\newline
Department of Mathematics, \newline
Integral University, Kursi Road, \newline
Lucknow-226026. \newline
Email : gazala.math@gmail.com\newline

\noindent
Mohd. Bilal\newline
Department of Mathematical Sciences, \newline
Umm Ul Qura University,\newline
Makkah, Saudi Arabia.\newline
Email: mohd7bilal@gmail.com


\begin{thebibliography}{99}

%\bibitem{AP} Alegre, P., Slant submanifolds of Lorentzian Sasakian and Para Sasakian
%manifolds, Taiwanese J. Math., 17 (2013), 897-910. 

%\bibitem{BR1}  Barbosa, E.,  Ribeiro, E., On conformal solutions of the Yamabe flow, Arch. Der Math., 101 (2013), 79-89.

%\bibitem{BR2}  Barros, A., Ribeiro, E., Some characterizations for compact almost Ricci solitons. Proc. Am. Math. Soc. 2012, 140, %1033-1040.

\bibitem{BAM1} Blaga, A. M., Solitons and geometrical structure in a perfect fluid spacetime, Rocky Mt. J. Math. (2020).

\bibitem{BAM2}  Blaga, A. M.,  Some geometrical aspects of Einstein, Ricci and Yamabe solitons, J. Geom. Symmetry Phys., 52 (2019), 17-26.

 \bibitem{CM1} Catino, G. and Mazzieri, L., Gradient Einstein solitons, Nonlinear Anal., 132(2016), 66-94.

\bibitem{CM2}  Catino, G., Cremaschi, L., Djadli, Z., Mantegazza, C. and Mazzieri, L., The Ricci Bourguignon flow, Pacific J. Math., 28(2017), 337-370.

\bibitem{CV} Chidananda, S., and Venkatesha, V., Yamabe soliton and Riemann soliton on Lorentzian para-Sasakian manifold, Commun. Korean Math. Soc., https://doi.org/10.4134/CKMS.c200365.

 %\bibitem{DC} Deshmukh, S. and Chen, B. Y.,  A note on Yamabe solitons,  Balkan Journal of Geometry and Its Applications, 23 (1)%(2018), 37-43.

%\bibitem{GDAB} Ganguly, D., Dey, S., Ali, A. and  Bhattacharyya, A., Conformal Ricci soliton and quasi-Yamabe soliton on %generalized Sasakian space form, Journal of Geometry and Physics, 169(2021), 104339.

\bibitem{GC} G\"{u}ler, S. and Crasmareanu, M., Ricci-Yamabe maps for Riemannian flows and their volume variation and volume entropy, Turk. J. Math., 43 (2019), 2631-2641.

\bibitem{HRS1} Hamilton,  R. S., Lectures on Geometric Flows (Unpublished manuscript, 1989).

\bibitem{HRS2} Hamilton, R. S., The Ricci Flow on Surfaces, Mathematics and General Relativity (Santa Cruz, CA, 1986), Contemp. Math., A.M.S., 71 (1988), 237-262.

%\bibitem{HC} Haseeb, A. and Chaubey, S. K., Lorentzian para-Sasakian manifolds and $*-$Ricci solitons, Kragujevac Journal of %Mathematics, 48(2) (2024), 167-179.

%\bibitem{HCM} Haseeb, A., Chaubey, S. K. and Khan, M. A., Riemannian $3$-manifolds and Ricci-Yamabe solitons, Int. J. Geom. %Met.  Mod. Phy. https://doi.org/10.1142/S0219887823500159.

\bibitem{HD} Haseeb, A. and De, U. C., $\eta$-Ricci solitons in $\epsilon$-Kenmotsu manifolds, J. Geom. 110, 34 (2019).

\bibitem{HA} Haseeb, A. and Almusawa, H., Some results on Lorentzian para-Kenmotsu manifolds admitting $\eta$-Ricci solitons, Palestine Journal of Mathematics, 11(2)(2022), 205-213.

\bibitem{HP1} Haseeb, A. and Prasad, R., Certain results on Lorentzian para-Kenmotsu manifolds, Bol. Soc. Parana. Mat., 39(3) (2021), 201-220.

\bibitem{HP2} Haseeb, A. and Prasad, R., Some results on Lorentzian para-Kenmotsu manifolds, Bull. Transilvania Univ. of Brasov, 13(62) (2020), no. 1, 185-198.

%\bibitem{HPCV} Haseeb, A., Prasad, R., Chaubey, S. K and Vanli, A. T., A note on $*$-conformal and gradient $*$-conformal %$\eta$-Ricci solitons in $\alpha$-cosymplectic manifolds, Honam Mathematical J., 44(2) (2022), 231-243.

\bibitem{HPM} Haseeb, A., Prasad, R. and Mofarreh, F., Sasakian manifolds admitting $*$-$\eta$-Ricci-Yamabe
solitons, Advances in Mathematical Physics, Vol. 2022, Article ID 5718736, 7 pages. doi: https://
doi.org/10.1155/2022/5718736


\bibitem{LHM} Li, Y., Haseeb, A. and  Ali, M., $LP$-Kenmotsu manifolds admitting $\eta$-Ricci solitons and spacetime, Journal of Mathematics, 2022, Article ID 6605127, 10 pages.

\bibitem{LH} Lone, M. A. and Harry, I. F.,  Ricci Solitons on Lorentz-Sasakian space forms, Journal of Geometry and Physics,
104547, doi: https://doi.org/10.1016/j.geomphys.2022.104547.

%\bibitem{IE} Ikawa, T. and  Erdogan, M., Sasakian manifold with Lorentzian metric, Kyungpook Math. J., 35(3) (1996), 517-526.

%\bibitem{OB} O'Neill, B., Semi-Riemannian Geometry with Applications to Relativity. Academic Press, New York (1983).


\bibitem{PCP} Pankaj, Chaubey, S. K and Prasad, R., Three dimensional  Lorentzian para-Kenmotsu manifolds and Yamabe soliton, Honam Mathematical J., 43(4) (2021), 613-626.

\bibitem{SK} Singh, J. P. and Khatri, M., On Ricci-Yamabe soliton and geometrical structure in a perfect fluid spacetime, Afr. Mat., 32(2021), 1645-1656. 
 
%\bibitem{SD} Suh, Y. J. and De, U. C., Yamabe Solitons and Ricci Solitons on almost co-Kahler manifolds, Canad. Math. Bull., 62 %(2019), 653-661.

\bibitem{VK} Venkatesha, Kumara, H. A., Ricci soliton and geometrical structure in a perfect fluid spacetime with torse-forming vector field, Afr. Mat. 30 (2019), 725-736 

\bibitem{YHI} Yoldas, H. I., On Kenmotsu manifolds admitting $\eta$-Ricci-Yamabe solitons,Int. J. Geom. Met.  Mod. Phy., 18(12) (2021), 2150189.

\bibitem{YK1} Yano, K., Integral Formulas in Riemannian geometry, Pure and Applied Mathematics, Vol. I, Marcel Dekker, New York, 1970.

\bibitem{YK2} Yano, K. and Kon, M., Structures on manifolds, World Scientific, (1984).

%\bibitem{ZL1} Zhang P., Li Y., Roy S. and De S., Geometry of $\alpha$-Cosymplectic Metric as $*$-Conformal
%$\eta$-Ricci-Yamabe Solitons Admitting Quarter-Symmetric Metric Connection, Symmetry, 13(11)(2021), 2189.

%\bibitem{ZL2}  Zhang, P., Li, Y., Roy, S., Dey, S. and Bhattacharyya, A., Geometrical Structure in a Perfect
%Fluid Spacetime with Conformal Ricci-Yamabe Soliton, Symmetry, 14(2022), 594.

\end{thebibliography}
\end{document}